%Paper: alg-geom/9505004
%From: shimada@mpim-bonn.mpg.de (Ichiro Shimada)
%Date: Thu, 4 May 1995 11:17:34 +0200
%Date (revised): Thu, 4 May 1995 12:01:28 +0200
 
\magnification=\magstep1
%macros
\font\BBF=msbm10
\def\Bbf#1{\hbox{\BBF #1}}
\font\SCR=cmsy10
\def\Scr#1{\hbox{\SCR #1}}

\def\Ker{\hbox {\rm Ker }}
\def\C{\Bbf C}
\def\G{\Bbf G}
\def\P{\Bbf P}

\def\Z{\Bbf Z}
\def\A{\Bbf A}
\def\O{\Scr O}
\def\qed
{\hskip 10pt \hbox{
\vrule height 7.5pt depth -0.1pt \vrule height 7.53pt depth -7.2pt width 7.3pt
\hskip -7.5pt \vrule height 0.3pt depth 0pt width 7.6pt \vrule height 7.5pt
depth -0.1pt
}}

\def\tl{\tilde}
\def\wt{\widetilde}

\def\sprime{\sp{\prime}}

\def\inv{\sp{-1}}

\def\im{\hbox{\rm im }}
\def\Im{\hbox{\rm Im }}

\def\Ker{\hbox{\rm Ker }}

\def\Coker{\hbox{\rm Coker }}
\def\dim{\hbox{\rm dim }}

\def\mod{ \hskip 5pt \hbox{\rm mod }}

\def\Sing{\hbox{\rm Sing\hskip 2pt}}
\def\aa{\alpha}
\def\bb{\beta}
\def\dd{\delta}
\def\ee{\epsilon}
\def\gg{\gamma}
\def\ll{\lambda}
\def\mm{\mu}
\def\nn{\nu}
\def\rr{\rho}
\def\ss{\sigma}

\def\zz{\zeta}

\def\xx{\xi}
\def\tt{\tau}
\def\ee{\epsilon}
\def\vee{\varepsilon}

\def\kk{\kappa}

\def\ii{\iota}
\def\vph{\varphi}

\def\SS{\Sigma}

\def\medn{\medskip\noindent}

\def\parn{\par\noindent}

\def\lra{\longrightarrow}

\def\lmt{\longmapsto}

\def\hs{\hskip 5pt}

\def\maprightsp#1{\smash{\mathop{\longrightarrow}\limits\sp{#1}}}

\def\varmaprightsp#1#2{\smash{\mathop{\hbox to #1 {\rightarrowfill}}
\limits\sp{#2}}}
\def\varmaprightsb#1#2{\smash{\mathop{\hbox to #1 {\rightarrowfill}}
\limits\sb{#2}}}

\def\mapdown#1{\Big\downarrow
\rlap{$\vcenter{\hbox{$\scriptstyle#1$}}$}}

\def\hookdownarrow
{{}\sp\cap \hskip -3.703pt \lower 2pt \hbox{$\downarrow$}}
\def\hookuparrow
{\lower 2pt\hbox{${}\sb{\cup}$}\hskip -3.7pt \lower -1pt\hbox{$\uparrow$}}
\def\veq{\|}
\def\downsim{\downarrow\hskip -2.3pt \wr}

\def\setbar{\hs ; \hs}
\def\set#1#2{\{\hs{#1}\setbar{#2}\}}
\def\seth#1#2{\{\hs{#1}\setbar{\hbox{#2}}\}}
\def\sethd#1#2#3
{\Bigl\{ \hs{#1}\setbar
{\matrix{\hbox{#2} \hfill\cr\hbox{#3} \hfill }} \Bigr\} }
\def\setht#1#2#3#4
{\biggl\{ \hs{#1}\setbar
{\matrix{\hbox{#2} \hfill\cr\hbox{#3} \hfill\cr\hbox{#4} \hfill}} 
\biggr\} }
\def\locus#1{\{ #1 \}}

\def\diagram#1{
\def\normalbaselines{\baselineskip20pt\lineskip3pt\lineskiplimit3pt}
\matrix{#1}}

\def\varvarmatrix#1#2#3#4
{\def\normalbaselines{\baselineskip#1\lineskip3pt\lineskiplimit3pt}
\vbox{\vskip #2 \hbox{\hfill$\matrix{#3}$\hfill}\vskip #3}}
\def\mapdiagram#1{
\def\normalbaselines{\baselineskip20pt\lineskip3pt\lineskiplimit3pt}
\matrix{#1}}
\def\varmatrix#1#2
{\def\normalbaselines{\baselineskip#1\lineskip3pt\lineskiplimit3pt}
\matrix{#2}}
\def\bigcases#1
{\biggl\{\, \vcenter{\normalbaselines{\mathsurround=0pt}
\ialign{$##\hfil$&\quad##\hfil\crcr#1\crcr}}\biggr.}
\def\st{\subset}
\def\sm{\setminus}

\def\pione{\pi\sb1}
\def\ratmap{\cdots\to}
\def\and{\hbox{and}}
\def\where{\hbox{where}}
\def\for{\hbox{for}}
\def\exact{\hbox{(exact)}}
\def\bfclaim#1{\medskip\noindent{\bf #1 }\hs }
\def\hsp#1{\phantom{#1}}
\def\shs#1{\hskip #1pt}
\def\afcite{\hsp{a}}
\def\mult{\sp{\times}}
\def\Pt{\P\sp 2}
\def\hc{ \xx\sb 0 : \xx\sb 1 : \xx\sb 2 }
\def\ac{ \xx\sb 0 , \xx\sb 1 , \xx\sb 2 }
\def\nt{ f\sb 1, \dots, f\sb n }
\def\nr{\sb{nr}}
\def\Sym{\hbox{Sym\hskip 2pt}}
\def\L{\Scr{L}}
\def\resU{|\sb{U} }
\def\resL{|\sb{L\sm\{ \infty \}}}
\def\resLC{|\sb{L\sm ( C\cup \{\infty \})}}
\def\AL{L\sm \{\infty \}}
\def\ALC{L\sm ( C\cup \{\infty \})}
\def\Gm{\G\sb m}
\def\trgp{\{ 1 \}}
\def\nz{\sp{\times}}
\def\ArtalBartolo{[1]}
\def\BruceGiblin{[2]}
\def\Dimca{[3]}
\def\FultonLazarsfeld{[4]}
\def\Libgober{[5]}
\def\Nemethi{[6]}
\def\OkaN{[7]}
\def\OkaS{[8]}
\def\OkaPre{[9]}
\def\ShimadaF{[10]}
\def\Shimadafinite{[11]}
\def\Tokunaga{[12]}
\def\Turpin{[13]}
\def\ZariskiE{[14]}
\def\ZariskiP{[15]}
\def\ZariskiT{[16]}
\def\DimcaPart#1{[3, #1]}
\def\FultonLazarsfeldPart#1{[4, #1]}
\def\ShimadaFPart#1{[10, #1]}
\def\link{\ell}
%
%
%text
%
%
\centerline{\bf A weighted version of Zariski's hyperplane section theorem }
\centerline{\bf and fundamental groups of complements of  plane curves}
\bigskip
\centerline{Ichiro Shimada}
\bigskip\noindent
{\bf \S 0. Introduction}
\medskip
In this paper, we formulate and prove 
a weighted homogeneous version of Zariski's hyperplane section theorem
on the fundamental groups of the complements of hypersurfaces
in a complex projective space,
and apply it to the study of 
$\pione (\Pt\sm C)$, where  $C\st \Pt$ is a projective plane curve.
The main application is to prove a comparison theorem as follows.
Let $\vph : \Pt \to \Pt$ be the composition of the Veronese embedding
$\Pt \hookrightarrow \P\sp N$
and the restriction of a general projection $\P\sp N \cdot\cdot \to \Pt$.
Our  comparison theorem
enables us to calculate 
$\pione (\Pt\sm\vph\inv(C))$ from $\pione (\Pt\sm C)$.
In \ZariskiE\afcite and \ZariskiT,
Zariski studied some projective plane curves 
with interesting properties.
An example is sextic curves with $6$ cusps.
Zariski showed that the fundamental group of the complement
depends on the placement of the $6$ cusps.
Another example is the $3$-cuspidal quartic curve,
whose complement has a non-abelian and finite fundamental group.
This curve is the only known example with this property.
Using the  comparison theorem, we derive infinite series 
of curves with these interesting properties from the classical examples
of Zariski.
As another application,
we shall discuss a  relation between 
$\pione(\Pt \sm C)$
and $\pione (\Pt \sm (C\cup L\sb{\infty}))$, 
where $L\sb{\infty}$ is a  line intersecting $C$ transversely.
\medskip
Let $x\sb 1$, \dots, $x\sb n$ be variables with weights
$$
\deg x\sb 1 = d\sb 1, \quad \dots \quad,  \deg x\sb n = d\sb n,
$$
and let $F(x\sb 1, \dots, x\sb n)$ be 
a non-zero weighted homogeneous polynomial
of total degree $d > 0$.
Suppose that $n \ge 2$ and $d\sb i > 0$ for $i= 1, \dots , n$.
In the affine space 
$\A\sp n$ with affine coordinates $(x\sb 1, \dots, x\sb n)$,
the equation $F=0$ defines a hypersurface $\SS \st \A\sp n$.
We let the multiplicative group $\G\sb m$ of non-zero complex numbers
act on $\A\sp n$ by
$$
\ll\cdot(x\sb 1, \dots, x\sb n) \quad = \quad
(\ll\sp{d\sb 1}x\sb 1, \dots , \ll\sp{d\sb n} x\sb n)
\qquad \for \quad \ll \in \G\sb m .
\eqno{(0.1)}
$$
This action leaves the complement $\A\sp n \sm \SS$ invariant,
because $F$ is weighted homogeneous.
Thus we have a natural homomorphism 
$$
\pione (\G\sb m) \quad \lra \quad \pione (\A\sp n \sm \SS).
\eqno{(0.2)}
$$
Note that the image of (0.2) is contained in the center of 
$\pione (\A\sp n \sm \SS)$,
so that we have the cokernel of this homomorphism.
We fix a projective plane $\Pt$ with homogeneous coordinates $(\hc )$.
Let $f\sb i$ be homogeneous polynomials of degree $d\sb i$ in $(\ac )$
for $i=1, \dots ,n $.
We denote by $f$ the $n$-tuple $ (\nt )$.
The polynomial 
$$F\sb f (\ac) := F(f\sb 1 (\ac), \dots, f\sb n (\ac) )$$
in the variables $(\ac )$ is then
homogeneous of degree $d$,
and $F\sb f (\ac )=0$ defines a subscheme  of $\Pt$,
which we shall denote by $C\sb f$.
Unless $F\sb f$ is constantly zero on $\Pt$,
$C\sb f$ is a projective plane curve of degree $d$.
\bfclaim{Theorem 1.}
{\sl 
Suppose that $\SS$ is reduced.
If $\nt$ are general,
then $\pione (\Pt \sm C\sb f)$ is isomorphic to the cokernel
of the natural homomorphism (0.2).
}
\medskip
Consider the case when $d\sb 1 = \cdots = d\sb n =1 $ and $n > 2$.
Then 
$F=0$ defines a hypersurface $\overline \SS$ 
in the projective space $\P\sp{n-1}$
with homogeneous coordinates $(x\sb 1: \dots : x\sb n)$.
In this case, the action (0.1) of $\Gm$ on $\A\sp n \sm \SS$
is fixed-point free,
and the quotient space $(\A\sp n \sm \SS) / \Gm$ is 
isomorphic to $\P\sp{n-1}\sm \overline \SS$.
Hence the cokernel of the homomorphism (0.2) is isomorphic to 
$\pione (\P\sp{n-1}\sm \overline \SS)$.
On the other hand, if $\nt$ are general linear forms,
then $C\sb f$ is the pull-back of $\overline \SS$
by the linear embedding $\ii\sb f : \Pt \hookrightarrow \P\sp{n-1}$
defined by $\ii\sb f\sp{*} {\shs 3} x\sb i = f\sb i$.
Hence Theorem 1 is nothing but the classical hyperplane section theorem of
Zariski \ZariskiP\afcite in this case.
This justifies us in calling Theorem 1
{\sl a weighted Zariski's hyperplane section theorem}.
\medskip
Let $f\sb 0, f\sb 1 , f\sb 2 $
be general homogeneous polynomials of degree $k$ in $(\ac)$.
We consider the branched covering $\vph :\Pt \to \Pt$ of degree $k\sp 2$
defined by
$$
(\hc) \quad\lmt\quad ( \hs f\sb 0 (\ac) \hs : 
\hs f\sb 1 (\ac ) \hs : \hs f\sb 2 (\ac) \hs ).
$$ 
Suppose that we are given 
a reduced projective plane curve $C \st \Pt$ of degree $d$.
Using Theorem 1,
we shall show that,
when $f\sb 0, f\sb 1, f\sb 2$ are general with respect to $C$,
the fundamental group $\pione (\Pt \sm \vph\inv (C))$
can be computed from $\pione (\Pt \sm C)$ in a simple way.
\par
We define the {\sl  linking number map}
$$
\link \quad : \quad \pione (\Pt \sm C) \quad \lra \quad \Z/ (d)
$$
as follows. 
It is well known that $H\sb 1 (\Pt \sm C , \Z)$
is naturally isomorphic to the cokernel of 
the direct sum of the restriction maps
$$
H\sp 2 (\Pt, \Z)
\quad \maprightsp{r\sb 1 \oplus \cdots \oplus r\sb s}\quad 
H\sp 2 (C\sb 1 , \Z) \oplus \cdots \oplus H\sp 2 (C\sb s , \Z),
\eqno{(0.3)}
$$
where $C\sb 1 , \dots , C\sb s$
are the irreducible components of $C$ (cf. \FultonLazarsfeldPart{\S 8} ).
Let $e\sb i \in H\sp 2 (C\sb i , \Z) \cong \Z$
be the positive generator, which is 
the Poincar\'e dual of a point on $C\sb i$,
and let $[L] \in H\sp 2 (\Pt, \Z) \cong \Z$
be the positive generator,
which is the Poincar\'e dual of a line.
Since $C$ is reduced, we have 
$r\sb i ([L]) = (\deg C\sb i ) \cdot e\sb i$.
Because $d=\deg C\sb 1  + \cdots + \deg C\sb s$,
the maps $e\sb i \mapsto 1 \mod d$ induces a well-defined homomorphism
from the cokernel of (0.3) to $\Z/(d)$; i.e.,
$H\sb 1 (\Pt \sm C, \Z) \to \Z/ (d)$.
The map $\link$ is then defined to be the composition of 
this homomorphism with the Hurwicz map
$\pione (\Pt \sm C) \to H\sb 1 (\Pt \sm C, \Z)$.
\par
We can define $\link$ in the following different way.
Let $[\aa] \in \pione (\Pt \sm C)$
be an element represented by a loop $\aa : S\sp 1 \to \Pt \sm C$.
Since $\pione (\Pt) =\trgp$,
there exists a continuous map $\bb : D \to \Pt$
from $D:=\set{z\in \C}{ |z| \le 1}$ to $\Pt$
such that $\partial \bb = \aa$.
By a small deformation of $\bb$ homotopically relative to the boundary,
we may assume that the image of $\bb$ intersects the curve $C$ transversely.
Let $\Im\bb\cdot C$ be the intersection number.
Since the intersection number of any element in $H\sb 2 (\Pt, \Z)$ with 
$[C] \in H\sb 2 (\Pt, \Z)$
is divisible by $d$,
the number   $\Im\bb\cdot C$ modulo $d$
is independent of the choice of $\bb$.
We define $\link([\aa])$ to be $\Im\bb\cdot C \mod d \in \Z/ (d)$.
From this definition,
we can call $\link([\aa])$
the linking number of the loop $\aa$ around $C$.
It is easy to see that this definition coincides with the previous definition.
\par
We define the {\sl extended linking number map} by
$$
\mapdiagram{
\tl \link &\quad:\quad& \Z\times \pione (\Pt \sm C) & \lra & \Z / (d) \cr
&& (\nu , [\aa]) & \mapsto & \link([\aa])-\nu \mod d .\cr
}
$$
\bfclaim{Theorem 2.}
{\sl
Let $C \st \Pt$  be  a reduced curve of degree $d$.
Suppose that $f\sb 0$, $f\sb 1$ and  $f\sb 2 $ are 
general homogeneous polynomials
of degree $k$ in $(\ac)$.
Then
$\pione (\Pt \sm \vph\inv(C))$ is isomorphic to 
$\Ker \tl \link / ((kd) \times\trgp )$.
}
\bfclaim{Corollary 1.}
{\sl
Suppose the same assumptions as in Theorem 2.
\par
(1)
The fundamental group  $\pione (\Pt \sm \vph\inv (C))$ is 
a central extension of $\pione (\Pt \sm C)$ by $\Z/ (k)$.
In particular,
if $\pione (\Pt \sm C)$ is finite,
then $\pione (\Pt \sm \vph\inv (C))$
is also finite of order
$k$ times the order of $\pione (\Pt \sm C)$.
\par
(2)
The fundamental group  $\pione (\Pt \sm \vph\inv (C))$ is abelian 
if and only if
so is $\pione (\Pt \sm C)$.
}
\medskip
Remark that,
when $f\sb 0$, $f\sb 1$ and $f\sb 2$ are general,
the morphism $\vph : \Pt \to \Pt$ is \'etale
over a Zariski open neighborhood of the singular points of $C$.
Hence, for example,
if the singular locus
$\Sing C$ consists of $\dd$ nodes and $\kk$ cusps,
then the singular locus of $\vph\inv (C)$ 
consists of  $k\sp 2 \dd$ nodes and $k\sp 2 \kk$ cusps.
\medn
{\sl Examples.}
\medskip
(1) {\sl Zariski pairs.}
\hs
A couple of reduced projective plane curves $C\sb 1$ and $C\sb 2$
is said to make a {\sl Zariski pair}
if they satisfy the following conditions \ArtalBartolo;
(i)
$\deg C\sb 1 = \deg C\sb 2$, 
(ii)
there exist tubular neighborhoods $T(C\sb i) \st \Pt $ of $C\sb i$ for $i=1,2$
such that $(T(C\sb 1), C\sb 1)$ and $(T(C\sb 2 ), C\sb 2)$
are diffeomorphic, and  (iii)
the pairs $(\Pt, C\sb 1)$ and $(\Pt, C\sb 2)$ are {\sl not} homeomorphic.
That is,  
the singularities of $C\sb 1 $ and $C\sb 2$ are topologically equivalent, 
but the embeddings $C\sb 1 \hookrightarrow \Pt$
and $C\sb 2 \hookrightarrow \Pt$ are not topologically equivalent.
\smallskip
The first example of Zariski pair was discovered and studied by Zariski.
In \ZariskiE\afcite  and \ZariskiT, he showed that
there exist projective plane curves $C\sb 1 $ and $C\sb 2$ of degree $6$
such that $\Sing C\sb 1$ consists of $6$ cusps lying on a conic,
while $\Sing C\sb 2$ consists of $6$ cusps {\sl not} on any conic, and that
$\pione (\Pt \sm C\sb 1)$ is isomorphic to the free product 
$\Z/(2) * \Z/(3)$ of cyclic groups of order $2$ and $3$,
while $\pione (\Pt \sm C\sb 2)$ is cyclic of order $6$.
Thus $C\sb 1$ and $C\sb 2$ make a Zariski pair.
(See also \OkaS.)
After this example,
only few Zariski pairs have been constructed (cf. \ArtalBartolo, \Tokunaga).
\par
Let $C\sb 1 (k)$ and $C\sb 2 (k)$
be the pull-backs of the sextic curves
$C\sb 1 $ and $C\sb 2$ above by the covering $\vph: \Pt \to \Pt$
in Theorem 2.
Both of them are of degree $6k$, and 
each of their singular loci consists of $6 k\sp 2$ cusps.
Combining Corollary 1 (2) with Zariski's result,
we see that $\pione (\Pt \sm C\sb 1 (k))$ is non-abelian, 
while $\pione (\Pt \sm C\sb 2 (k))$
is abelian.
Thus we obtain an infinite series of Zariski pairs 
$C\sb 1 (k)$ and $C\sb 2 (k)$.
\medskip
(2)
{\sl Pull-backs of the three cuspidal quartic.}
\hs
Let $C\sb 0 \st \Pt$ be a curve of degree $4$ 
with $3$ cusps and no other singularities;
for example, the curve defined by
$$
x\sp 2 y\sp 2 + y\sp 2 z\sp 2 + z\sp 2 x\sp 2 - 2 xyz(x+y+z) = 0.
$$
In fact, it is known
that any three cuspidal quartic curve $C\sb 0$ is projectively isomorphic
to the curve defined by this  equation \BruceGiblin.
This curve was discovered and studied by Zariski in \ZariskiE.
(See also  \DimcaPart{Chapter 4, \S 4}.)
The remarkable property of $C\sb 0$ is that $\pione (\Pt \sm C\sb 0)$
is isomorphic to the binary dihedral group of order $12$.
Other than this three cuspidal quartic,
there have been no examples of projective plane curves $C\st \Pt$
such that $\pione (\Pt \sm C)$ is non-abelian and finite.
\par
Let $C\sb 0 (k)$ be the pull-back of $C\sb 0$ by $\vph :\Pt \to \Pt$ in Theorem 2.
Then $C\sb 0 (k)$ is of degree $4k$ and $\Sing C\sb 0 (k)$
consists of $3k\sp 2$ cusps.
By Corollary 1,
$\pione (\Pt \sm C\sb 0 (k))$
is non-abelian and finite of order $12k$.
\par 
In the forthcoming paper \Shimadafinite,
we will construct
other examples of such  curves.
\medskip
(3) {\sl The fundamental group $\pione (\Pt \sm C\sb{p,q,k})$.}
\hs
Let $p$ and $q$ be positive integers prime to each other,
and let $f$ and $g$ be general homogeneous polynomials in $(\ac)$
of degree $pk$ and $qk$, respectively, with $k \ge 1$.
We define   $C\sb{p,q,k}$ 
to be  the projective plane curve of degree $pqk$ 
defined by $f\sp q + g\sp p = 0$ (cf. \Libgober).
The fundamental group of $U:=\A\sp 2 \sm \{ x\sp q + y \sp p = 0 \}$
is well-known (cf. \DimcaPart{Chapter 4, \S 2}),
and the homomorphism from $\pione (\Gm)$ to $\pione (U)$
induced by the action of
$\Gm$ on $U$ with  weights $(pk, qk)$ on variables $(x,y)$
can be easily described.
Then, from Theorem 1, we see that $\pione (\Pt \sm C\sb{p,q,k})$ is 
isomorphic to
$$
\langle  \hs  a, b, c \hs  | 
\hs a\sp q = b \sp p = c, \hs c \sp k = 1 \hs  \rangle.
$$
Some parts of this fact have been already proved;  
by Zariski \ZariskiE\afcite when $p=2$, $q=3$ and $k=1$,
by Turpin \Turpin \afcite when $p=2$, $q=3$ and  $k>1$,
and by Oka \OkaN \afcite and N\'emethi \Nemethi\afcite 
when $p, q $ are arbitrary and $k=1$.
\medskip
As another  corollary of the proof of Theorem 2,
we will show the following:
\bfclaim{Corollary 2.}
{\sl Let $C$ be a reduced
projective plane curve of degree $d$, and let $L\sb{\infty}$ be a general line.
Then the fundamental group $\pione (\Pt\sm (C \cup L\sb{\infty}))$
of the affine part of $\Pt \sm C$ 
is isomorphic to
$\Ker \tl \link$.
In particular,
$\pione (\Pt\sm (C\cup L\sb{\infty}))$ is abelian
if and only if so is $\pione (\Pt \sm C)$.
}
\medskip
Oka raised the following problem:
\medn
{\sl Question.}
Let $C\sb 1$ and $C\sb 2$ be  projective plane curves
such that $\pione (\Pt \sm C\sb 1 ) $ and $\pione (\Pt \sm C\sb 2 )$ are 
isomorphic.
Let $L$ be a general line.
Then are the fundamental groups 
$\pione (\Pt \sm (C\sb 1 \cup L)) $ and $\pione (\Pt \sm (C \sb 2 \cup L) )$
isomorphic?
\medn
Suppose that $C\sb 1 $ and $C\sb 2 $ are reduced and
irreducible.
If the isomorphism 
$H\sb 1 (\Pt \sm C\sb 1, \Z ) \cong H\sb 1 (\Pt \sm C\sb 2, \Z )$, 
induced from 
the given isomorphism 
$\pione (\Pt \sm C\sb 1 ) \cong \pione (\Pt \sm C\sb 2 )$ 
via the Hurwicz homomorphism, maps the positive generator of 
$H\sb 1 (\Pt \sm C\sb 1, \Z ) $
to the positive generator of $H\sb 1 (\Pt \sm C\sb 2, \Z )$,
then the answer is affirmative because of Corollary 2.
\medskip
The main tool of the proof of Theorem 1 is \ShimadaFPart{Theorem 1}.
A similar idea was used to calculate $\pione (\Pt \sm C\sb{p,q, 1})$ in
\ShimadaFPart{\S 4}.
Theorem 2 is proved  by applying Theorem 1 
to the case when $n= 3$, $d\sb 1 = d\sb 2 = d\sb 3 =k$,
and $F$ is the homogeneous polynomial defining   $ C\st \Pt$.
\medskip
Independently from us,
Oka \OkaPre\afcite also composed examples of projective plane curves as 
our examples (1) and (2) above.
As our method,
his construction makes use of the covering of the plane,
but not of the projective plane as ours,
but an affine part of it.
The curves he constructed have singularity 
of different types from our examples.
The method of proof is also quite different.
\bfclaim{Acknowledgment.}
The results of this paper have been 
obtained in an effort to answer various problems discussed at
the workshop on ``Fundamental groups and branched covering"
held at Tokyo Institute of Technology on December 1994.
The author thanks to the Professor M.\ Oka  
for inviting me to this workshop,
and for stimulating discussions.
\bigskip\noindent
{\bf \S 1. Proof of the weighted Zariski's hyperplane section theorem}
\medskip
We consider $\Pt$ as the space of  $1$-dimensional linear subspaces 
in a $3$-dimensional linear space $V$ with linear coordinates $(\ac)$.
Let
$$
A :=  H\sp 0 (\Pt, \O(d\sb 1)) \times \cdots \times H\sp 0 (\Pt, \O(d\sb n))
$$
be the space of all $n$-tuples $f= (\nt)$.
Then we have a natural morphism
$$
\mapdiagram{
 \Psi &\quad : \quad & V\times A & \lra & \A\sp n \cr
&& (\hs (\ac), \hs (\nt) \hs ) & \lmt &
(\hs f\sb 1 (\ac), \hs  \dots \hs , \hs  f\sb n (\ac) \hs ). \cr
}
$$
Let $\wt W \st V \times A$ be the pull-back of $\SS \st \A\sp n$ by $ \Psi$.
Since $F$ is weighted homogeneous,
we can make  $\Gm$ act on the complement 
$(V\times A) \sm \wt W$
by
$$
\ll\cdot ((\ac), (\nt)) =  
(\hs (\ll\xi\sb0, \ll\xi\sb1, \ll\xi\sb 2), \hs  (\nt)\hs )
\quad\where\quad \ll \in \Gm. \eqno{(1.1)}
$$
The morphism
$$
 \Phi \quad : \quad (V\times A ) \sm \wt W \quad \lra \quad \A\sp n \sm \SS,
$$
which is the restriction of $ \Psi$,
is equivariant under the actions (0.1) and (1.1) of $\Gm$ on each side.
Let $ W \st \Pt \times A$
be the divisor defined by
$$
F(f\sb 1 (\ac), \dots, f\sb n (\ac)) = 0,
$$
which is the universal family of the subschemes 
$\{ C\sb f  \}\sb {f \in A}$
of $\Pt$ parameterized by $A$.
The action (1.1) of $\Gm$ on $(V\times A ) \sm \wt W$ is fixed-point free, 
and the quotient space is nothing but the complement $(\Pt \times A) \sm W$.
Hence we get a commutative diagram
$$
\diagram{
\pione (\Gm) & \lra & \pione ((V\times A) \sm \wt W) & \lra &
 \pione ((\Pt \times A) \sm W) & \lra & \trgp &\quad\exact \cr
\veq & & \phantom{\Phi\sb{*}} \downarrow \Phi\sb{*} & & && \cr
\pione (\Gm) & \lra &\pione (\A\sp n \sm \SS) & \lra &
 \hbox{the cokernel of (0.2)} & \lra & \trgp &\quad\exact. \cr
}
$$
On the other hand,
the complement 
$\Pt \sm C\sb f$ is the fiber over 
$f \in A$ of the second projection $(\Pt \times A) \sm W \to A$.
Hence Theorem 1 follows from the following two claims.
\bfclaim{Claim 1.}
{\sl
The homomorphism $ \Phi\sb{*} : 
\pione ((V\times A)\sm\wt W) \to \pione (\A\sp n \sm \SS)$
is an isomorphism.
}
\bfclaim{Claim 2.}
{\sl
The inclusion $\Pt\sm C\sb f \hookrightarrow (\Pt \times A) \sm W$
induces an isomorphism on the fundamental groups when $f \in A$ is 
general enough.
}
\medskip
{\it Proof of Claim 1.}
\hs
Note that $\{ o\sb V \} \times A \st V\times A$
is contained in $\wt W$,
where $o\sb V $ is the origin of $V$.
Using the homotopy exact sequence,
we can verify Claim 1 by proving that
the restriction of $ \Psi$ to
$(V\sm \{o\sb V \}) \times A$
gives $(V\sm \{o\sb V \})\times A$
a structure of the locally trivial fiber space
over $\A\sp n$
with simply connected fibers.
Note that the restriction of $\Psi$  to $\{v\} \times A$ is a surjective
affine linear map from 
$\{ v \} \times A \cong A$ 
to $\A\sp n$ for any point $v \in V\sm\{ o\sb V \}$.
For arbitrary points $u \in \A\sp n$ and $v \in V\sm \{ o\sb V \}$,
the intersection $ \Psi \inv (u) \cap (\{ v \} \times A)$
is then isomorphic, via the second projection,
to an affine subspace of $A$ of codimension $n$.
Thus every fiber of 
$ \Psi |\sb{ ( V\sm \{ o\sb V\} ) \times A } : 
( V\sm \{o\sb V \}) \times A \to \A\sp n$
is a fiber bundle over $V\sm \{ o \sb V \}$ with fibers 
isomorphic to an affine space of dimension $\dim A - n$.
\qed
\medskip
{\it Proof of Claim 2.}
\hs
Let $\Xi\nr \st A$ be the locus of all $f \in A$
such that $C\sb f$ is {\sl not} a reduced curve;
that is, $f=(\nt)$ belongs to $\Xi\nr$
if and only if $F\sb f (\ac)$ is constantly zero on $\Pt$,
or $F\sb f (\ac) = 0$ defines a non-reduced curve.
By \ShimadaFPart{Theorem 1},
the proof of Claim 2 is reduced to the verification of the following:
\bfclaim{Claim 2'.}
{\sl
The locus $\Xi\nr \st A$ is of codimension $\ge 2$.
}
\medskip
Before proving this claim, we need two preparations.
\medskip
Let $\SS\sb 1 $, \dots  , $\SS\sb M$ be the irreducible components of 
$\SS \st \A\sp n$,
and let $(\Sing \SS)\sb 1 $, \dots  , $(\Sing \SS)\sb{M\sprime}$ 
be the irreducible components of $\Sing \SS$ of dimension $n-2$.
We give them the reduced structure as subschemes of $\A\sp n$.
For $f \in A$,
let $ \psi \sb f : V\sm \{ o\sb V \} \to \A\sp n$
be the restriction of $ \Psi : V\times A \to \A\sp n$
to the subvariety 
$(V\sm \{ o\sb V \} ) \times \{ f \} \cong V\sm \{ o\sb V \}$.
\bfclaim{Sub-claim.}
{\sl
Let $f \in A$ be generally chosen.
Then, for $i = 1, \dots, M$ and $j = 1, \dots , M\sprime $,
the following conditions ($i$) and ($j\sprime$) hold.
\par
($i$)
\hs
There exists a point $P\sb i \in V\sm \{ o\sb V \}$
such that $ \psi \sb f (P\sb i ) $ is contained in $\SS\sb i \sm (\Sing \SS)$,
and the pull-back of $\SS$ by $ \psi\sb f$ is non-singular and 
of codimension $1$
locally around $P\sb i$.
\par
($j\sprime $)
\hs
There exists a point $P\sb j \sprime \in V \sm \{ o\sb V \}$
such that $ \psi\sb f (P\sb j \sprime)$
is contained in $(\Sing \SS)\sb j$,
and the pull-back of $(\Sing \SS)\sb j$ by $ \psi \sb f $ is 
non-singular and of codimension $2$ locally around $P\sb j \sprime$.
}
\medn
{\it Proof.}\hs
The conditions $(i)$ and $(j\sprime)$
are open  for the choice of $f \in A$.
Therefore, it is enough to show that,
for each condition,
there exists at least one $f \in A$ which satisfies it.
This can be verified from the following fact.
Let $R$ be the point $(1,0,0)$ on $V \sm \{ o\sb V \}$.
For an arbitrary point $Q = (c\sb 1, \dots, c\sb n) \in \A\sp n$
and arbitrary two tangent vectors 
$$
\def\tanvec#1{{\partial\over{\partial x\sb{#1}}}}
{\bf u } = u\sb 1 \tanvec 1 + \cdots + u\sb n 
\tanvec n \in T\sb Q (\A \sp n) , 
\quad\and\quad
{\bf v } = v\sb 1 \tanvec 1 + \cdots + v\sb n \tanvec n \in T\sb Q (\A \sp n)
$$
of $\A\sp n$ at $Q$,
there exists an element $f \in A$
such that $ \psi \sb f  (R) = Q$,
and the image of the tangential map
$ \psi\sb{f*} : T\sb R ( V\sm \{ o\sb V \} ) \to T\sb Q (\A\sp n)$
contains both of  ${\bf u}$  and  ${\bf v}$.
Indeed, 
if the coefficients of $\xi\sb 0 \sp{d\sb i}$, 
$\xi\sb 0 \sp{d\sb i -1}\xi\sb 1$ and 
$\xi\sb 0 \sp{d\sb i -1} \xi\sb 2 $ in $f\sb i$
are $c\sb i$, $u\sb i$ and $v\sb i$, respectively,
then $f = (\nt)$ satisfies the required condition. \qed
\medskip
Next we consider the action of the general linear transformation group
$GL (V)$ of $V$ on the space $A$.
Note that $H\sp 0 (\Pt, \O (1))$ can be identified with
the dual space $V\sp * := \hbox{Hom } (V, \C)$.
We let $GL(V)$ act on $V\sp *$ from left
by
$$
(g(l))(v) = l (g\inv (v)) \quad\for\quad g \in GL(V),  \hs l \in V\sp * \hs
\and \hs  v \in V.
$$
This action can be extended in the natural way to the action on
$A = \Sym\sp{d\sb 1} V\sp * \times \cdots \times \Sym\sp{d\sb n} V\sp *$.
Thus $GL(V)$ acts on $(V\sm \{ o\sb V \}) \times A$, and $\Pt \times A$.
By definition,
we have
$$
 \Psi (g( P), g(f)) =  \Psi ( P, f) \eqno{(1.2)}
$$
for every $ P \in V\sm \{ o\sb V \} $, $f \in A$ and $g \in GL (V)$.
In particular, 
the divisor $W \st \Pt \times A$ is invariant under the action of $GL(V)$,
and $C\sb{g(f)} = g(C\sb f)$ for every $g \in GL(V)$ and $f\in A$.
It follows that $\Xi\nr \st A$ is also invariant under the action of $GL(V)$.
\medskip
Now we start the proof of Claim 2'.
For a line $L \st \Pt$,
we put
$$
\Xi (L) := \seth{f \in A}{ $L\cap C\sb f$ does {\it not} 
consist of distinct $d$ points}.
$$
The locus $\Xi\nr$ is contained in $\Xi (L)$
for every line $L$.
We fix a line $L\sb 0 \st \Pt$,
and let $\Xi (L\sb 0)\sb1$, \dots , $\Xi(L\sb 0)\sb N$
be the irreducible components of $\Xi (L\sb 0)$.
Note that $\Xi (g(L\sb 0)) = g(\Xi (L\sb 0))$
for every $g \in GL(V)$. 
Since the subgroup 
$G(L\sb 0) \st GL (V)$
of all $g \in GL(V)$ 
which leave $L\sb 0$ invariant
is connected,
the action of $g \in G(L\sb 0)$ on $\Xi (L\sb 0)$
does not interchange irreducible components of $\Xi (L\sb 0)$;
that is, each irreducible component is invariant
under the action of $G(L\sb 0)$.
Thus,
for an arbitrary line $L$,
we have a natural numbering $\Xi (L)\sb1$, \dots , $\Xi(L)\sb N$
of the irreducible components of $\Xi(L)$
such that $g(\Xi (L)\sb i) = \Xi (L\sb 0)\sb i$
for every $g \in GL(V)$ with $g(L) = L\sb 0$.
We shall show the following:
\medn
{\bf Claim 2".}
\hs
{\sl
For each $i = 1, \dots, N$,
$\Xi (L\sb 0)\sb i \st A$
is of codimension $\ge 1$,
and $\Xi (L\sb 0)\sb i \sm \Xi\nr $
is non-empty.
}
\medn
Since $\Xi\nr$ is a Zariski closed subset of 
$\Xi(L\sb 0)$, this will prove Claim 2'.
Roughly speaking, the idea of the proof of Claim 2"  is to show that,
if $f\sp 0 \in A$ is general,
then $f\sp 0 \not \in \Xi\nr$,
while there exist lines $L\sb i $ $(i= 1, \dots, N)$
such that $f\sp 0 \in \Xi (L\sb i) \sb i$.
Suppose that
these are proved.
Let $g\sb i \in GL(V)$ be a linear transformation  such that 
$g\sb i (L\sb i) = L\sb 0$.
Then, we have 
$g\sb i ( \Xi (L\sb i)\sb i \sm \Xi\nr) = \Xi(L\sb 0)\sb i \sm \Xi\nr$,
and this contains an element $g\sb i (f\sp 0)$,
so that 
$\Xi (L\sb 0)\sb i \sm \Xi\nr$ is non-empty.
\medskip
To carry out this idea, we have to investigate the irreducible components
of $\Xi(L\sb 0)$ more closely.
Let $P$ be a point on a line $L$.
We put
$$
\Xi(L, P) := \sethd{f\in A}{the restriction of $F\sb f (\ac)
 \in H\sp0(\Pt, \O(d))$
to $L$ has }{a zero of order $\ge 2$ at $P$} .
$$
As before, we have 
$$
g(\Xi(L, P)) = \Xi(g(L), g(P)) \eqno{(1.3)}
$$
for every $g \in GL(V)$.
We also have
$$
\Xi (L\sb 0) = \bigcup\sb{P\in L\sb 0} \Xi (L\sb 0, P) 
= \bigcup \sb {g(L\sb 0) = L\sb 0} g (\Xi (L\sb 0, P\sb 0)),
$$
where $P\sb 0 \in L\sb 0$ is a fixed point.
From these two,
by the same argument as above,
it follows that the set of irreducible components of $\Xi (L\sb 0, P\sb 0)$
corresponds bijectively and in a natural way to that of $\Xi (L\sb 0)$.
Let $\Xi (L\sb 0, P\sb 0)\sb 1$, \dots , $\Xi (L\sb 0, P\sb 0)\sb N$
be the numbering of the irreducible components of $\Xi(L\sb 0, P\sb 0)$
according to that
of $\Xi (L\sb 0)$; that is, 
we have $\Xi (L\sb 0, P\sb 0)\sb i \st \Xi (L\sb 0 )\sb i$ for 
$i = 1, \dots , N$.
Since $\dim \Xi(L\sb 0) \le \dim \Xi(L\sb 0, P\sb 0) +1$,
it is enough to show that,
for $i=1, \dots, N$,
the locus $\Xi (L\sb 0, P\sb 0) \sb i \st A$ is of codimension $\ge 2$
and $\Xi(L\sb 0, P\sb 0)\sb i \sm \Xi\nr$ is non-empty.
\medskip
Let $L\sb 0$ be defined by $\xi \sb 2 = 0$,
and  $P\sb 0$  the point $(1:0:0)$ on $L\sb 0$.
Let $t= \xi\sb 1 / \xi \sb 0$ be the affine coordinate on $L\sb 0$
with the origin $P\sb 0$. We put
$$
G\sb f (t) := F\sb f (1, t, 0) =
F(f\sb 1 (1, t , 0), \dots, f\sb n (1, t, 0)).
$$
Then $f$ is contained in $\Xi (L\sb 0, P\sb 0)$
if and only if
$$
G\sb f (0) = {{dG\sb f}\over{dt }}(0) = 0 \eqno{(1.4)}
$$
holds.
Let $a\sb i$ be the coefficient of $\xi\sb 0 \sp{d\sb i}$ in $f\sb i$,
and $b\sb i$ the coefficient of 
$\xi\sb 0\sp{d\sb i - 1} \xi\sb 1 $ in $f\sb i$.
Then $(a\sb 1, \dots, a\sb n, b\sb 1, \dots, b\sb n)$
form a subset of a linear coordinate system of $A$.
Let $\A\sp n \times \A\sp n$
be the affine space with affine coordinates 
$(a\sb 1,\dots, a\sb n, b\sb 1, \dots, b\sb n)$.
There exists a natural projection $q : A \to \A\sp n \times \A\sp n$.
The condition (1.4)
can be written
as follows in terms of  $(a\sb 1,\dots, a\sb n, b\sb 1, \dots, b\sb n)$.
$$
\def\pF#1{{{\partial F}\over{\partial x\sb{#1}}}}
\eqalign
{
F(a\sb 1, \dots, a\sb n) = 0 , \quad \and \hskip 80pt\cr
\pF{1}(a\sb 1, \dots, a\sb n)\cdot b\sb 1 + \cdots 
+ \pF{n}(a\sb 1, \dots, a\sb n)\cdot b\sb n = 0.
}
$$
We consider $\A\sp n \times \A\sp n$
with the first projection
$\A\sp n \times \A\sp n \to \A\sp n$
as the tangent bundle $T\A\sp n$ of the first factor $\A\sp n$
in which the zero section is given by
$P\mapsto (P,O)$
({\it not} by the diagonal $P\mapsto (P,P)$)
where $O$ is the origin 
$(b\sb 1, \dots, b\sb n) = (0, \dots, 0)$
of the second factor.
Under this identification,
the  equations above define in $\A\sp n \times \A\sp n$
the space 
$$
T\SS  := 
\sethd{(P, Q) \in \A\sp n \times \A\sp n}{$Q$ is contained in 
the Zariski tangent }{space $T\sb P (\SS) \st T\sb P (\A\sp n)$ 
of $\SS$ at $P$} .
$$
Thus $\Xi(L\sb 0, P\sb 0)$ can be identified with 
$q\inv (T\SS)$,
which is isomorphic to the product of $T\SS$
with an affine space of dimension $\dim A - 2n$.
Since $\SS$ is reduced by the assumption,
$T\SS$ is of dimension $2n-2$,
and hence $\Xi(L\sb 0, P\sb 0) \st A$ is of codimension $2$.
Thus $\Xi (L\sb 0) \st A$ is of codimension $\ge 1$,
and hence so is $\Xi\nr \st A$.
\medskip
Let $p : T\SS \to \SS$ be the  projection 
induced by the natural projection
$T \A\sp n \cong \A\sp n \times \A\sp n \to \A\sp n$.
The fiber of $p$ over $P \in \SS$ is a linear space of dimension $n-1$ 
(resp. $n$) if $\SS$ is non-singular (resp. singular) at $P$.
Thus the irreducible components of $T\SS$ are listed
as follows; the closure of $p\inv (\SS\sb i \sm (\Sing \SS) )$ 
for $i = 1, \dots , M$,
and $p\inv ((\Sing \SS)\sb j)$ for $j = 1, \dots, M\sprime$.
Hence we get $N= M+ M\sprime$, and by changing the numbering, we have
$$
\varmatrix{18pt}{
\Xi(L\sb 0, P\sb 0)\sb i\hsp{aaa} & =& 
q\inv(\hbox{ the closure of $p\inv (\SS\sb i \sm (\Sing \SS) )$ }) 
& \hs\for\hs &i = 1, \dots , M, &\and \cr
\Xi(L\sb 0, P\sb 0)\sb{M+j} & =& q\inv(p\inv ((\Sing \SS)\sb j))\hfill& 
\hs\for\hs& j = 1, \dots , M\sprime&. \cr
} 
$$
Let $\wt P\sb 0 \sp{(1)}$ be the point $(1,0,0) \in V\sm \{ o\sb V \}$,
which is located over $P\sb 0 \in \Pt$.
Note that $f \in \Xi (L\sb 0, P\sb 0)$
implies $P\sb 0 \in C\sb f$,
and hence $\psi\sb f (\wt P\sb 0 \sp{(1)})\in \SS$.
We identify two affine spaces 
$\A\sp n$ with coordinates $(a\sb 1, \dots, a\sb n)$ and 
$\A\sp n$ with coordinates $(x\sb 1, \dots, x\sb n)$
by putting
$a \sb i = x \sb i$
for $i=1, \dots, n$.
Then, recalling the definitions of the morphisms
$ \psi\sb f$, $q$ and $p$,
we can easily check that
$$
\psi \sb f (\wt P\sb 0 \sp{(1)} ) = p(q(f)) \quad 
\hbox {for every}\quad f \in \Xi (L\sb 0, P\sb 0).
$$
Note that $\psi\sb f : V\sm \{o\sb V \} \to \A \sp n$ is 
equivariant under the actions of $\Gm$
with weights $(1,1,1)$ on the left-hand side and 
with weights $(d\sb 1, \dots , d\sb n)$ on the right-hand side.
Moreover, each of 
$\SS\sb i \st \A \sp n$ and $(\Sing \SS)\sb j \st \A\sp n $ are 
invariant under the action of $\Gm$,
because $\Gm$ is connected.
Therefore, given an element $f \in \Xi (L\sb 0, P\sb 0)$,
we can tell which irreducible components of 
$\Xi (L\sb 0, P\sb 0)$ this $f$ belongs to
by looking at
the point $\psi \sb f (\wt P\sb 0 ) \in \SS$
where $\wt P\sb 0 \in V\sm \{o\sb V \}$
is an arbitrary point positioned over $P\sb 0 = (1:0:0) \in \Pt$.
\medn
{\sl Criterion.}
\hs
If $\psi \sb f (\wt P\sb 0 ) \in \SS$
is a non-singular point of $\SS\sb i$,
then $f$ belongs to $\Xi (L\sb 0, P\sb 0)\sb i$.
If
$\psi \sb f (\wt P\sb 0 ) \in \SS$
is contained in $(\Sing \SS )\sb j$,
then $f$ belongs to $\Xi (L\sb 0, P\sb 0)\sb{M+j}$.
\medskip
Recall that we have already proved that $\Xi\nr \st A$ is of
codimension $\ge 1$.
By Sub-claim,
there exists $f\sp 0 \in A\sm \Xi\nr$
such that the image of $\psi\sb{f\sp 0}$
intersects $\SS\sb i \sm (\Sing \SS)$ for
$i=1, \dots, M$,
and $(\Sing \SS)\sb j$ for $j=1, \dots, M\sprime$.
Let $\wt P\sb i$ and $\wt P\sb j\sprime$
be points on $V\sm \{o\sb V \}$
such that $\psi\sb{f\sp 0} (\wt P\sb i) \in \SS\sb i \sm (\Sing \SS)$
and $\psi\sb{f\sp 0} (\wt P\sb j\sprime) \in (\Sing \SS)\sb j$.
We denote by $P\sb i $ and $P\sb j\sprime$ the points on $\Pt$
corresponding to $\wt P\sb i$ and $\wt P\sb j\sprime$, respectively.
It is obvious that these points are on $C\sb{f\sp 0}$.
There exists a line $L\sb i \st \Pt$
which intersects
$C\sb{f\sp0}$ at $P\sb i$
with multiplicity $\ge 2$.
This means that $f\sp 0 \in \Xi (L\sb i, P\sb i)$.
Let $g\sb i \in GL(V)$ be a linear transformation 
such that $g\sb i (P\sb i) = P\sb 0 $
and $g\sb i (L\sb i) = L\sb 0$.
Then, by (1.3),  we have $g\sb i (f\sp 0) \in \Xi (L\sb 0, P\sb 0)$.
By (1.2),
we also have
$$
 \psi \sb{g\sb i (f\sp 0)} (g\sb i (\wt P\sb i)) =  
\psi\sb{f\sp 0} (\wt P\sb i) \in \SS\sb i \sm (\Sing \SS).
$$
Since $g\sb i (\wt P\sb i) \in V\sm \{o\sb V \}$ is 
a point over $ g\sb i ( P\sb i) = P\sb 0$,
we see that $g\sb i (f\sp 0)$ belongs to $\Xi (L\sb 0, P\sb 0) \sb i$
by the above criterion.
From $f\sp 0 \not \in \Xi\nr$,
we have $g\sb i (f\sp 0) \not \in \Xi\nr$.
Thus $\Xi (L\sb 0, P\sb 0) \sb i \sm \Xi\nr$
is non-empty for $ i = 1, \dots, M$.
Similarly,
there exists a line $L\sb j\sprime \st \Pt$
which intersects
$C\sb{f\sp0}$ at  $P\sb j\sprime$ with multiplicity $\ge 2$,
which means  $f\sp 0 \in \Xi (L\sb j\sprime, P\sb j\sprime)$.
We choose $g\sb j\sprime \in GL(V)$ 
such that $g\sb j\sprime (P\sb j\sprime) = P\sb 0 $
and $g\sb j\sprime (L\sb j\sprime) = L\sb 0$.
Then (1.3) implies $g\sb j\sprime (f\sp 0) \in \Xi (L\sb 0, P\sb 0)$.
By (1.2),
we  have
$ \psi \sb{g\sb j\sprime (f\sp 0)} (g\sb j\sprime (\wt P\sb j\sprime)) 
=  \psi\sb{f\sp 0} (\wt P\sb j\sprime) \in (\Sing \SS)\sb j$.
This means   $g\sb j\sprime (f\sp 0) \in \Xi (L\sb 0, P\sb 0) \sb{ M+j}$
by the criterion above.
Thus $\Xi (L\sb 0, P\sb 0) \sb {M+j} \sm \Xi\nr$ 
is non-empty for $ j = 1, \dots, M\sprime$.
\qed
\bigskip\noindent
{\bf \S 2. Proof of the comparison theorem}
\bigskip
We shall prove Theorem 2 and Corollaries in this section.
In fact, most parts of Corollaries can be proved directly
without using Theorem 2.
\medskip
Suppose that $C\st \Pt$
is defined by the homogeneous equation $F(\ac)=0$.
Let $\SS\st \A\sp 3$ be the hypersurface defined by $F=0$
in the affine space $\A\sp 3$
with affine coordinates $(\ac)$.
For a positive integer $\nu$,
we denote by 
$$
\rr\sb{\nu} \quad : \quad \pione (\Gm) \quad \lra \quad 
\pione (\A\sp 3 \sm \SS)
$$
the homomorphism 
induced by the action of $\Gm$ on $\A\sp 3 \sm \SS$ with weights 
$(\nu, \nu, \nu)$;
$$
\ll\cdot (\ac) =
(\ll\sp{\nu}\xi\sb 0, \ll\sp{\nu}\xi\sb 1, \ll\sp{\nu}\xi\sb 2)
\quad\where\quad \ll\in\Gm.
\eqno{(2.1)}
$$
The image of $\rr\sb{\nu}$ is 
contained in the center of $\pione(\A\sp 3 \sm \SS)$.
It is obvious that $\pione (\Pt\sm C)$ is isomorphic to
the cokernel of $\rr\sb 1$.
Since $C$ is reduced by the assumption,
$\SS$ is also reduced.
Now, applying Theorem 1 with $n=3$ and $d\sb 1 =d\sb 2 =d\sb 3 =k$,
we see that $\pione (\Pt \sm \vph \inv (C))$
is isomorphic to the cokernel of $\rr\sb k$,
because $\vph\inv(C)$ is defined by $F(f\sb 0, f\sb 1, f\sb 2) =0$.
Let 
$$
\ss\sb{\nu} \quad  : \quad  \pione (\Gm) \quad\lra\quad  \pione (\Gm)
$$
be the homomorphism induced from the morphism
$\Gm\to \Gm$ given by $\ll \mapsto \ll\sp{\nu}$.
Then we get a commutative diagram
$$
\diagram{
\pione(\Gm) & \maprightsp{\rr\sb k} & \pione (\A\sp 3 \sm \SS ) & \lra & 
\pione (\Pt\sm \vph\inv (C)) & \lra & \trgp & 
\quad \exact\phantom{ .}  \cr
\mapdown{\ss \sb k}&& \veq &&&&&&\cr
\pione(\Gm) & \maprightsp{\rr\sb 1} & \pione (\A\sp 3 \sm \SS ) & \lra & 
\pione (\Pt\sm C) & \lra & \trgp & 
\quad \exact .\cr
}
\eqno{(2.2)}
$$
From this diagram, we can naturally derive an exact sequence
$$
\Coker \ss\sb k \hs \maprightsp{\dd} \hs  \pione (\Pt\sm \vph\inv (C)) 
\hs \lra \hs \pione (\Pt\sm C)
\hs \lra \hs \trgp \quad \exact. \eqno{(2.3)}
$$
The cokernel of $\ss\sb k$ is isomorphic to $\Z / (k)$.
Hence
(2.3) proves the implications
$$\vbox{
\halign{
#\hfill&\quad # \quad&#\hfill \cr
$\pione (\Pt \sm C)$ is non-abelian 
&$\Longrightarrow$&
$\pione (\Pt \sm \vph\inv(C))$ is non-abelian, \qquad \and \cr
$\pione (\Pt \sm C)$ is finite
&$\Longleftrightarrow$&
$\pione (\Pt \sm \vph\inv(C))$ is finite. \cr
}}\eqno{(2.4)}
$$
Thus half of Corollary 1 is already proved.
\smallskip
Since the image of $\rr\sb 1 $ is contained in
the center
of $\pione (\A\sp 3 \sm \SS)$,
the image of $\dd$ in (2.3) is also contained in the center
of $\pione (\Pt \sm \vph\inv(C))$.
We shall show that $\dd$ is injective,
and the central extension (2.3)
is the same as the one described in Theorem 2.
For this purpose, we investigate the bottom  exact sequence of (2.2) 
more closely.
\medskip
From now on,
we denote by $U$ the complement $\Pt \sm C$.
For a positive integer $m$,
let $\L\sb{m}\to \Pt$ be the line bundle corresponding to the invertible sheaf
$\O(-m)$ on $\Pt$,
and let $\L\nz\sb{m} \st \L\sb{m}$ denote the complement of the zero section.
Let $O$ be the origin of $\A\sp 3$, and 
let $\A\sp 3 \sm\{ O \} \to \Pt$ be the natural projection;
that is,
the quotient map by the $\Gm$-action (2.1) with $\nu = 1$.
There exists
an isomorphism between $\L\sb 1 \nz$ and $\A\sp 3 \sm\{O\}$
over $\Pt$.
We fix such an isomorphism.
From this, we get an isomorphism 
$$
\diagram{
\L\sb1\nz\resU &&\cong&& \A\sp 3 \sm \SS \cr
&\searrow&&\swarrow& \cr
&&U&&\cr
}
$$
over $U$,
where $\L\sb 1\nz\resU$ is 
the restriction of $\L\sb 1\nz \to \Pt $ to $U\st \Pt$,
and $\A\sp 3 \sm \SS \to U$ is the quotient map by the $\Gm$-action (2.1)
with $\nu = 1$.
Hence we have an isomorphism between the homotopy exact sequences;
$$
\diagram{
\lra&\pione(\Gm) & \maprightsp{\rr\sb 1} & \pione (\A\sp 3 \sm \SS ) & \lra & 
\pione (U) & \lra & \trgp & 
\quad \exact\phantom{ .}  \cr
&\downsim&&\downsim&& \veq &&&\cr
\lra &\pione(\C\nz) & \maprightsp{\ii\sb 1} & \pione (\L\sb1\nz\resU) & \lra & 
\pione (U) & \lra & \trgp & 
\quad \exact .  \cr
}
\eqno{(2.5)}
$$
Here we denote the fiber of $\L\sb 1\nz \to \Pt$
over the base point of $U$ by $\C\nz$.
\medskip
Since the line bundle $\L\sb 1 \sp{\otimes d}$ is isomorphic to $\L\sb d$, 
there exists a morphism  
$$
\mm \quad : \quad  \L\sb1\nz \quad \lra \quad \L\sb d\nz \eqno{(2.6)}
$$
over $\Pt$
which induces $\zz\sb 1\mapsto\zz\sb d =\zz\sb 1\sp d$ on
the  fibers $\C\nz$ of $\L\nz \sb 1$ and $\L\nz\sb d$ over 
the base point of $U$
with appropriate fiber coordinates $\zz\sb1$ and $\zz\sb d$.
From this morphism,
we get a homomorphism 
between the homotopy exact sequences 
of $\L\sb1\nz \resU \to U$ and $\L\sb d\nz\resU \to U$;
$$
\diagram{
\lra&\pione(\C\nz) & \maprightsp{\ii\sb 1} & \pione ( \L\sb1\nz\resU) & 
\maprightsp{\tt\sb 1} & 
\pione (U) & \lra & \trgp & 
\quad \exact\phantom{ ,}  \cr
&\mapdown{\mm\sb{\#}}&&\mapdown{\mm\sb{*}}&& \veq &&&\cr
\lra &\pione(\C\nz) & \maprightsp{\ii\sb d} & \pione (\L\sb d\nz\resU) & 
\maprightsp{\tt\sb d} & 
\pione (U) & \lra & \trgp & 
\quad \exact ,  \cr
}
\eqno{(2.7)}
$$
where $\mm\sb{\#}$ is the multiplication by $d$ on $\pione (\C\nz) \cong \Z$.
Note again that the images of $\ii\sb 1$ and $\ii\sb d$ are contained 
in the centers 
of $\pione (\L\sb 1\nz\resU) $
and $\pione (\L\sb d\nz\resU) $, respectively.
Since $\L\sb d \to \Pt$
is a line bundle corresponding to the invertible sheaf $\O(-d) \cong \O(-C)$,
there exists a meromorphic section
$s : \Pt \ratmap \L\sb d$
which has no zeros on $\Pt$ and is holomorphic outside $C$.
Restricting $s$ to $U$,
we get a holomorphic section 
$s : U\to \L\sb d \nz \resU$.
Hence the homotopy exact sequence of $\L\sb d\nz \resU \to U$
splits.
In particular,
the homomorphism $\ii\sb d$ in (2.7) is injective,
and we have an isomorphism
$$
(\ii\sb d , s\sb{*} )\quad : \quad \pione(\C\nz) \times \pione (U) 
\quad\maprightsp{\sim}\quad \pione (\L\sb d\nz \resU).
\eqno{(2.8)}
$$
Because $\mm\sb{\#}$ is also injective,
we see by diagram chasing that
$\ii\sb 1$ and $\mm\sb * $ are injective, too.
Since $\ii\sb1 $ is injective,
$\rr\sb 1$ in (2.5) is also injective 
and hence $\dd$ in (2.3) is injective.
Thus the proof of Corollary 1 (1) is completed.
\smallskip
We shall show that the image 
of the injective homomorphism $\mm\sb{*}$ in (2.7)
is a normal subgroup of $\pione (\L\sb d \nz \resU)$.
Let $[\aa]$ and $[\gg]$
be arbitrary elements of $\pione (\L\sb 1 \nz\resU)$ 
and $\pione (\L\sb d\nz\resU)$,
respectively,
and we put
$[\bb]:= [\gg]\inv\cdot \mm\sb{*} ([\aa])\cdot [\gg]$.
Let $[\gg\sprime ] \in \pione (\L\sb 1 \nz\resU)$
be an element such that $\tt\sb 1 ([\gg\sprime])= \tt\sb d ([\gg])$.
We put $[\aa\sprime]:=[\gg\sprime]\inv\cdot[\aa]\cdot[\gg\sprime]$ and 
$[\dd]:= [\gg]\inv \cdot\mm\sb{*}([\gg\sprime])$.
Then we have $\mm\sb{*} ([\aa\sprime ]) = [\dd]\inv[\bb][\dd]$.
Since $\tt\sb d ([\dd]) = 1$ by the definition,
$[\dd]$ is contained in the image of $\ii\sb d$,
which is in the center of $\pione (\L\sb d \nz\resU)$.
Hence we have $[\bb]=\mm\sb{*} ([\aa\sprime])$.
\smallskip
Now we can derive   the following commutative diagram from the diagram (2.7);
$$
\diagram{
&&\trgp&&\trgp&&&&&\cr
&&\mapdown{}&&\mapdown{}&&&&&\cr
\trgp &\lra&\pione(\C\nz) & \maprightsp{\ii\sb 1} & \pione ( \L\sb1\nz\resU) & 
\maprightsp{\tt\sb 1} & 
\pione (U) & \lra & \trgp & 
\quad \exact\phantom{ ,}   \cr
&&\mapdown{\mm\sb{\#}}&&\mapdown{\mm\sb{*}}&& \veq &&&\cr
\trgp&\lra &\pione(\C\nz) & \maprightsp{\ii\sb d} & \pione (\L\sb d\nz\resU) & 
\maprightsp{\tt\sb d} & 
\pione (U) & \lra & \trgp & 
\quad \exact ,   \cr
&&\mapdown{}&&\mapdown{\psi}&&&&&\cr
&&\Z/(d)&=&\Z/(d)&&&&&\cr
&&\mapdown{}&&\mapdown{}&&&&(2.9)&\cr
&&\trgp&&\trgp&&&&&\cr
&&\exact&&\exact.&&&&&\cr
}
$$
We write the definition of $\psi$ in this diagram explicitly.
For $[\gg] \in \pione (\L\sb d\nz\resU)$,
let $[\gg\sprime] \in \pione (\L\sb 1\nz\resU)$
be an element such that $\tt\sb 1 ([\gg\sprime]) = \tt\sb d ([\gg])$.
There exists a unique element $[\ee] \in \pione (\C\nz)$
such that $\ii\sb d ([\ee] ) = [\gg]\cdot \mm\sb * ([\gg\sprime] )\inv$.
We define
$$
\psi ([\gg]) \quad := \quad  [\ee] \mod \im \mm\sb{\#} \quad \in \quad \Z/(d).
$$
The independence of $\psi([\gg])$ on the choice of $[\gg\sprime]$
can be checked easily.
It is in
showing that $\psi$ is a group homomorphism that
we have to use the fact that
the image of $\ii\sb d$ is contained in the center.
Let $[\gg\sb 1]$ and $[\gg\sb 2]$ be two elements  in 
$\pione (\L\sb d\nz\resU)$.
Because 
$[\gg\sb 2]\cdot \mm\sb *([\gg\sb2 \sprime ] )\inv \in \im \ii\sb d$,
we have
$$
[\gg\sb1]\mm\sb * ([\gg\sb1\sprime ] )\inv \cdot [\gg\sb2] \mm\sb * 
([\gg\sb2\sprime ] )\inv
=[\gg\sb1][\gg\sb2]\cdot 
\mm\sb *([\gg\sb2\sprime ] )\inv\mm\sb *([\gg\sb1\sprime ] )\inv
=([\gg\sb1][\gg\sb2]) \mm\sb *([\gg\sb1\sprime ] [\gg\sb2\sprime ] )\inv.
$$
This implies $\psi([\gg\sb1])\psi([\gg\sb2]) = \psi([\gg\sb1][\gg\sb2])$. 
Now the surjectivity of $\psi$ and $\Ker \psi = 
\Im \mm\sb *$ can be checked immediately.
\medskip
The diagram (2.9)
shows that $\pione (\A\sp 3 \sm \SS)$,
which is isomorphic to
$\pione (\L\sb1\nz\resU)$
by  (2.5),
is isomorphic to the kernel of
$$
\psi \quad : \quad \pione (\L\sb d\nz\resU) \hs \cong\hskip -20pt\lower 14pt
\hbox{\rm by (2.8)}
\Z\times \pione (U) \quad \lra \quad \Z/(d).
$$
Because the image of $\ii\sb d \circ \mm\sb{\#}$ is equal 
to $(d)\times\trgp$ in
$\pione (\L\sb d\nz\resU) \cong \Z \times \pione (U)$,
the image of $\rr\sb 1 :\pione (\Gm) \to \pione (\A\sp 3 \sm \SS)$,
which corresponds to the image of
$\ii\sb 1 : \pione (\C\nz) \to \pione (\L\sb 1\nz\resU)$ 
via the isomorphism (2.5),
is given by
$$
(d) \times \trgp \quad \st \quad \Ker \psi \quad \cong \quad 
\pione (\A\sp 3 \sm \SS).
$$
By the diagram (2.2),
the image of $\rr\sb k : \pione (\Gm)  \to \pione (\A\sp 3 \sm \SS)$
is given by
$$
(kd) \times \trgp \quad \st \quad \Ker \psi \quad \cong \quad \pione 
(\A\sp 3 \sm \SS).
$$
Thus we obtain isomorphisms
$$
\pione (\Pt\sm C) \hs \cong \hs  \Ker\psi/ ((d) \times \trgp),
\quad\and\quad
\pione (\Pt\sm \vph\inv(C)) \hs \cong \hs \Ker\psi/ ((kd) \times \trgp).
\eqno{(2.10)}
$$
\par
Now we can complete the proof of Corollary 1.
The only remaining part is the implication
that, if $\pione (\Pt \sm C)$ is abelian,
then
so is $\pione (\Pt\sm \vph\inv (C))$.
If $\pione (\Pt\sm C) = \pione (U)$
is abelian,
then  $\Ker \psi \st \Z \times \pione (U)$ is also abelian,
and hence so is
$\pione (\Pt\sm\vph\inv(C)) \cong \Ker\psi /((kd) \times \trgp)$.
Note that we have also proved
that $\pione (\A\sp 3 \sm \SS)$
is abelian if and only if so is $\pione (\P\sp 2 \sm C)$.
\medskip
In order to prove Theorem 2,
it is enough to show that $\psi$
coincides with $-1$ times the extended linking number map $\tl \link$ 
defined in Introduction,
and this is equivalent to
show that the homomorphism
$$
\psi\circ s\sb{*} \quad : \quad \pione (U) \quad \lra \quad \Z/(d)
$$
coincides with $-\link$.
It is obvious that
$\psi\circ s\sb{*}$
factors
as follows;
$$
\pione (U) \quad\maprightsp{\eta}\quad  H\sb 1 (U,\Z) \cong 
\Coker (\shs 2 H\sp 2 (\Pt, \Z) \to 
\mathop\bigoplus\limits\sb{i=1}\sp{s} H\sp 2 (C\sb i, \Z) \shs 2) 
\quad \lra \quad\Z/(d),
$$
where $\eta$ is the Hurwicz map, and $C\sb1$, \dots, $C\sb s$ 
are the irreducible components of $C$.
By the definition of $\link$, it is enough to show that,
if $[\aa] \in \pione (U)$ is mapped to the $i$-th positive generator
$$
\ee\sb i := (0, \dots, 0,  e\sb i, 0, \dots , 0) \mod  H\sp 2 (\Pt, \Z) 
\quad \in \quad \Coker (\shs 2 H\sp 2 (\Pt, \Z) \to 
\mathop\bigoplus\limits\sb{i=1}\sp{s} H\sp 2 (C\sb i, \Z) \shs 2)
$$
by the Hurwicz map,
then $\psi\circ s\sb{*} ([\aa]) = -1 \mod d$.
Here $e\sb i \in H\sp 2 (C\sb i, \Z)$ is the Poincar\'e dual of 
a point on $C\sb i$.
\medskip
Before showing this,
we remark a property of the Hurwicz map.
The Hurwicz map does not depend on the choice of the base point.
Namely,
if we connect two base points
$b\sb 1$ and $b\sb 2$ in $U$ by a path,
we get an isomorphism $\pione (U, b\sb 1 ) \cong \pione (U, b\sb2)$,
which depends on the homotopy class of the connecting path.
But the  diagram
$$
\diagram{
\pione(U,b\sb 1 ) & & \cong & & \pione (U, b\sb 2) \cr
& \eta\searrow && \swarrow \eta & \cr
&& H\sb 1 (U, \Z)
}
$$
is always commutative,
whichever homotopy class of paths we may have chosen.
\medskip
Let $L\st \Pt$ be a general line,
and let $t$ be an affine coordinate on $L$
such that $t=\infty$ is not on the intersection $L\cap C$.
Note that $\L\sb {m} | \sb {L\sm \locus{\infty}}$ is a trivial 
bundle over $L\sm \locus{\infty}$ for any $m$.
Let $\zz\sb 1$ and $\zz\sb d$ be the fiber coordinates
of $\L\sb 1 \resL$ and $\L\sb d \resL$,
respectively, such that the morphism $\mm$ in (2.6) is given by 
$\zz\sb 1 \mapsto \zz\sb d = \zz\sb 1 \sp d$
over $L\sm \locus{\infty}$.
Suppose that $t=a\sb1, a\sb2, \dots , a\sb d$ are 
the intersection points of $L$ and $C$.
Then the section $s: U \to \L\sb d\nz\resU$,
restricted to $L\sm (C\cup\locus{\infty})$, is given by, for example,
$$
t \quad\lmt\quad (t, \zz\sb d ) = (t, f(t))
\quad\where \quad f(t) ={{1}\over{(t-a\sb1)\cdots(t-a\sb d)}}.
$$
Let us consider the section
$$
\mapdiagram{
s\sb 1 &\quad : \quad & \AL & \quad \lra \quad &\L\sb1\nz\resL \cr
&& t & \quad \lmt \quad & (t, \zz\sb 1 ) = (t, B) \hs , \cr 
} 
$$
where $B$ is a non-zero constant.
(Of course, this section $s\sb 1$ cannot extend over the whole $L$.)
We fix a base point $b$ of $U$.
By the remark above,
we may assume that $b$ is on $L$ and 
is close enough to an intersection point $t=a\sb{\nn}$ of $L$ 
and an irreducible component 
$C\sb i$ of $C$.
Let $t= t\sb b$ be the coordinate of the base point.
We choose the non-zero constant $B$ in such a way that
$f (t\sb b ) = B\sp d$;
that is,
the images of $s$ and $\mm\circ s\sb 1$ on 
the base point $b \in L$ coincide.
Let $\aa$ be the loop
$$
\mapdiagram{
\aa &\quad :\quad & [0, 1] & \quad \lra \quad & L \sm C \cr
&&\theta &\lmt & t = a\sb{\nn} + \vee e\sp{2\pi i \theta}, \cr
}
$$ 
where $\vee = t\sb b - a\sb{\nn}$.
It is easy to see that $\eta ([\aa])$ is the $i$-th positive generator 
$\ee\sb i$
corresponding to $C\sb i$.
We calculate $\psi\circ s\sb{*} ([\aa]) = \psi ([s\circ \aa])$.
It is obvious that
$[s\sb 1 \circ \aa] \in \pione (\L\sb 1 \nz \resU, s\sb 1 (b))$
satisfies $\tt\sb 1 ([s\sb1 \circ \aa ] ) = [\aa]$.
Hence $\psi ([s\circ \aa ]) $ is represented by
$[s\circ\aa]\cdot (\mm\sb * ([s\sb 1 \circ \aa ] )) \inv$,
which is in the image of $\ii\sb d$.
Recall that $\L\sb d\nz \to \Pt$ is trivial on $\ALC$;
that is, we have an isomorphism
$$
\L\sb d \nz \resLC \quad \cong \quad (\ALC ) \times \C\nz
\eqno{(2.11)}
$$
given by the coordinates $(t, \zz\sb d)$.
Hence we have
$$
\pione(\L\sb d \nz \resLC, s(b)) \quad \cong \quad
\pione  (\ALC , b ) \times \pione(\C\nz, B\sp d).
$$
In this direct product, we have
$$
[s\circ\aa ] = ([\aa], [\bb])\qquad\and \qquad 
\mm\sb* ([s\sb 1 \circ \aa ] ) = ([\aa], 0),
$$
where $\bb$ is the loop on $\C\mult$ obtained from the loop 
$s\circ\aa$ on $\L\sb d\nz\resLC$
by the second projection
$\L\sb d\nz\resLC \to \C\nz$
in (2.11), 
which can be written explicitly as follows;
$$
\theta \quad \lmt \quad f(a\sb{\nn} + \vee e\sp{2\pi i \theta} ) \in \C\mult .
$$
Their difference 
$[s\circ\aa]\cdot (\mm\sb * ([s\sb 1 \circ \aa ] )) \inv$
is then given by $[\bb] \in \pione (\C\nz)$.
Since $ | \vee |= | t\sb b - a\sb{\nn} |$ is small enough,
$\bb$ is homotopically equivalent in $\C\nz$
to the loop
$$
\theta\quad\lmt\quad B\sp d e\sp{-2\pi i \theta},
$$
which corresponds to $-1 \in \Z \cong \pione(\C\nz)$.
Thus $\psi\circ s\sb *([\aa]) = -1 \mod d$,
and Theorem 2 is proved.
\medskip
Now we shall prove Corollary 2.
We have already shown that $\pione (\A\sp 3 \sm \SS)$
is isomorphic to $\Ker \tl \link$,
and it is abelian if and only if so is $\pione (\Pt\sm C)$.
Hence
it is enough to show that,
when $L$ is a general line,
$\pione (\Pt\sm (C\cup L))$
is isomorphic to $\pione (\A\sp 3 \sm \SS)$.
We consider $\A\sp 3 $
as the complement of a hyperplane $H$ in $\P\sp 3$.
Let $\Pt \st \P\sp 3$ be a general plane.
Then $\Pt\cap (\SS\cup H)$ is isomorphic to $C\cup L$.
By the classical Zariski's hyperplane section theorem,
we have 
$\pione (\Pt\sm (C\cup L)) \cong \pione (\P\sp 3 \sm (\SS \cup H))
\cong \pione (\A\sp 3 \sm \SS)$.
\qed
\bigskip\noindent
\bigskip\noindent
\centerline{\bf References}
\medskip
\item{\ArtalBartolo}
Artal Bartolo, E.:
Sur les couples de Zariski.
J.\ Alg.\ Geom. {\bf 3}, 223 - 247 (1994)
\item{\BruceGiblin}
Bruce, J.W., Goblin, P.J.:
A stratification of the space of plane quartic curves.
Proc.\ London Math.\ Soc. {\bf 42}, 270 - 298 (1981)
\item{\Dimca}
Dimca, A.:
Singularities and Topology of Hypersurfaces,
Berlin Heidelberg New York: Springer 1992
\item{\FultonLazarsfeld}
Fulton, W., Lazarsfeld, R.:
Connectivity and its applications in
algebraic geometry.
(Lecture Notes in Math., vol. 862, pp. 29 - 92)
Berlin Heidelberg New York: Springer 1981
\item{\Libgober}
Libgober, A.:
Fundamental groups of the complements to plane singular curves.
Proc.\ Symp.\ in Pure Math.
{\bf 46}, 29 - 45 (1987)
\item{\Nemethi}
N\'emethi, A.:
On the fundamental group of the complement
of certain singular plane curves.
Math.\  Proc.\  Cambridge Philos.\  Soc. {\bf 102}, 453  -  457 (1987)
\item{\OkaN}
Oka, M.:
Some plane curves whose complements
have non - abelian fundamental groups.
Math.\  Ann. {\bf 218}, 55  -  65 (1975)
\item{\OkaS}
Oka, M.:
Symmetric plane curves with nodes and cusps.
J.\ Math.\ Soc.\ Japan {\bf 44}, 375  -  414 (1992)
\item{\OkaPre}
Oka, M.:
Two transformations of plane curves and their fundamental groups.
preprint
\item{\ShimadaF}
Shimada, I.:
Fundamental groups of open algebraic varieties.
to appear in Topology
\item{\Shimadafinite}
Shimada, I.: 
On projective plane curves whose complements
have finite non - abelian fundamental groups.
preprint
\item{\Tokunaga}
Tokunaga, H.:
A remark on Bartolo's paper.
preprint
\item{\Turpin}
Turpin, W.S.:
On the fundamental group of a certain class of plane curves.
Amer.\ J.\ Math. {\bf 59}, 529 - 577 (1937)
\item{\ZariskiE}
Zariski, O.:
On the problem of existence of 
algebraic functions of two variables
possessing a given branch curve.
Amer.\ J.\ Math. {\bf 51}, 305  -  328 (1929)
\item{\ZariskiP}
Zariski, O.:
A theorem on the Poincar\'e group of an algebraic hypersurface.
Ann.\ Math. {\bf 38}, 131  -  141 (1937)
\item{\ZariskiT}
Zariski, O.:
The topological discriminant group of a Riemann surface 
of genus $p$.
Amer.\ J.\ Math. {\bf 59}, 335 - 358 (1937)
\bigskip\noindent
Max-Planck-Institut f\"ur Mathematik 
\parn
Gottfried-Claren Str. 26
\parn
53225 Bonn, Germany
\parn
shimada@mpim-bonn.mpg.de
\end